\documentclass{amsart} 

\usepackage{amsmath,amsthm}     
\usepackage{amssymb}            
\usepackage{euscript}           
\usepackage{enumerate,calc}
\usepackage[matrix,arrow,curve,frame]{xy}    

\xymatrixcolsep{1.9pc}                          
\xymatrixrowsep{1.9pc}
\newdir{ >}{{}*!/-5pt/\dir{>}}                  

\newtheorem{thm}[subsection]{Theorem}
\newtheorem{defn}[subsection]{Definition}
\newtheorem{prop}[subsection]{Proposition}

\newtheorem{cor}[subsection]{Corollary}
\newtheorem{lemma}[subsection]{Lemma}

\theoremstyle{definition}  
\newtheorem{ex}[subsection]{Example}

\newtheorem{remark}[subsection]{Remark}

  {\end{list}}
  
\newcommand{\dfn}{\textbf} 

\newcommand{\mdfn}[1]{\dfn{\mathversion{bold}#1}} 

\newcommand{\tens}              {\otimes}               
\newcommand{\iso}               {\cong}  

\newcommand{\cat}{\EuScript}    
\newcommand{\cC}{{\cat C}}
\newcommand{\cD}{{\cat D}}

\newcommand{\cF}{{\cat F}}

\newcommand{\cP}{{\cat P}}

\newcommand{\cS}{{\cat S}}

\newcommand{\sSet}{s{\cat Set}}

\newcommand{\Ho}{\text{Ho}\,}
\newcommand{\ho}{\text{Ho}\,}
\newcommand{\et}{\text{et}}

\newcommand{\field}[1]  {\mathbb #1} 
\newcommand{\A}         {\field A}
\newcommand{\R}         {\field R}

\newcommand{\Z}         {\field Z}
\newcommand{\C}         {\field C}
\newcommand{\Q}         {\field Q}
\renewcommand{\P}       {\field P}
\newcommand{\F}         {\field F}

\DeclareMathOperator*{\colim}{colim}

\DeclareMathOperator{\spec}{Spec}
\DeclareMathOperator{\Spec}{Spec}
\DeclareMathOperator{\Hom}{Hom}
\DeclareMathOperator{\Ext}{Ext}
\DeclareMathOperator{\Aut}{Aut}

\DeclareMathOperator{\Map}{Map}

\DeclareMathOperator{\chara}{char}

\DeclareMathOperator{\Et}{Et}
\DeclareMathOperator{\SEt}{SEt}

\newcommand{\ra}{\rightarrow}                   
\newcommand{\lra}{\longrightarrow}              
\newcommand{\la}{\leftarrow}                    
\newcommand{\llra}[1]{\stackrel{#1}{\lra}}      

\newcommand{\inc}{\hookrightarrow}              

\newcommand{\blank}{-}                          
\newcommand{\id}{\textup{id}}                   

\newcommand{\rea}[1]{|{#1}|}             
\newcommand{\map}{\rightarrow}

\newcommand{\ceck}[1]{\Cech(#1)}         
\newcommand{\oceck}[1]{\Cech^{o}(#1)}    
\newcommand{\oreal}[1]{\rea{\oceck{U}}}  
\newcommand{\creal}[1]{\rea{\ceck{U}}}   

\newcommand{\Cech}{\check{C}}

\newcommand{\RP}{\R{\text{\sl P}}}
\newcommand{\CP}{\C{\text{\sl P}}}

\newcommand{\DQ}{DQ}

\numberwithin{equation}{subsection}

\newenvironment{myequation}
  {\addtocounter{subsection}{1}\begin{eqnarray}}
  {\end{eqnarray}$\!\!$}

\newcommand{\pro}{\textup{pro-}}
\newcommand{\pros}{\textup{pro}}

\newcommand{\Ecomp}{E^{\wedge}}
\newcommand{\SP}{\cS\cP}

\newcommand{\brF}{\overline{\F}}
\newcommand{\tmf}{tmf}

\begin{document}

\title{Etale homotopy and sums-of-squares formulas}
\author{Daniel Dugger}
\author{Daniel C. Isaksen}

\address{Department of Mathematics\\University of Oregon\\Eugene, OR 97403}
\address{Department of Mathematics\\Wayne State University\\Detroit, MI 48202}

\email{ddugger@math.uoregon.edu}
\email{isaksen@math.wayne.edu}

\begin{abstract}
This paper uses a relative of $BP$-cohomology to prove a theorem 
in characteristic $p$ algebra.
Specifically, we obtain some new
necessary conditions for the existence of sums-of-squares formulas
over fields of characteristic $p>2$.  These conditions
were previously known in characteristic zero by results of Davis.  Our
proof uses a generalized \'etale cohomology theory called \'etale
$BP2$.
\end{abstract}

\maketitle

\tableofcontents

\section{Introduction}
\label{sctn:intro}

Let $F$ be a field.  A \dfn{sums-of-squares formula} of type $[r,s,n]$
is an identity of the form
\[ (x_1^2+\cdots+x_r^2)(y_1^2+\cdots+y_s^2)=z_1^2+\cdots +z_n^2 
\]
where each $z_i$ is a bilinear expression in the $x$'s and $y$'s with
coefficients in $F$.  The identity occurs in the polynomial
ring $F[x_1,\ldots,x_r,y_1,\ldots,y_s]$.  A classical problem, raised
by Hurwitz in 1898 \cite{H1}, is to determine which triples $[r,s,n]$ can
admit a sums-of-squares formula (with the answer possibly depending on
the field $F$).  In this paper we show how \'etale homotopy theory and
classical algebraic topology combine to yield necessary conditions for
the existence of such a formula, applicable even to fields of
characteristic $p>2$.  Specifically, we use the generalized cohomology
theory $BP2$ to obtain some conditions previously unknown in the
positive characteristic case.

\medskip

Hurwitz introduced the sums-of-squares problem after he proved that a
formula of type $[n,n,n]$ can exist only if $n$ is $1$, $2$, $4$, or
$8$ \cite{H1}.  This is his classical theorem on the existence of
composition algebras.  Technically speaking, Hurwitz proved this
theorem only when $F=\R$, but his methods work for arbitrary $F$ not
of characteristic $2$.  (Note that the whole sums-of-squares problem
is trivial in characteristic $2$, where formulas of all possible types
exist).  Later, Hurwitz and Radon \cite{H2,R} completely analyzed the
case of formulas of type $[r,n,n]$.

When $F=\R$, topological methods have provided a series of necessary
conditions for a formula of type $[r,s,n]$ to exist.  The earliest papers in
this vein are by Hopf \cite{H} and Stiefel \cite{St},
who showed that if a formula of type $[r,s,n]$ exists
over $\R$ then there is a map of spaces
\[ \phi\colon \RP^{r-1}\times \RP^{s-1} \ra \RP^{n-1} 
\]
which restricts to linear inclusions on $\RP^{r-1}\times \{*\}$ and
$\{*\} \times \RP^{s-1}$ (for any choice of basepoints).  Such a $\phi$ is
called an \dfn{axial map} of projective spaces.  Hopf used
singular $\Z/2$-cohomology to obtain topological obstructions to the
existence of such a map.

Topologists have maintained an interest in the existence of axial maps
because of its relation to the immersion problem for $\RP^k$.  
Generalized cohomology theories supply more and more
stringent conditions on the existence of axial maps, developed in
papers such as \cite{As}, \cite{D1}, \cite{D2}, and \cite{BDM}.
The
paper \cite{A} also dealt with the immersion problem---and,
implicitly, the existence of sums-of-squares formulas---but not 
through the axial map perspective.

These topological methods give conditions on the existence of
sums-of-squares formulas over $\R$.  As explained in \cite[Sec.~2]{L},
these conditions also apply to any field of characteristic zero.  The
present paper continues the sequence \cite{DI1, DI2} in which we prove
that these results also apply to characteristic $p$ fields.  The
topological methods are replaced by methods of motivic homotopy theory.
The paper \cite{DI1}
used mod $2$ motivic cohomology, and \cite{DI2} used algebraic
$K$-theory.  The present paper uses \'etale $BP$-theory.

For additional background on the sums-of-squares problem over characteristic
zero fields, the reader may consult \cite{L} and \cite{Sh}.

\medskip

We now describe the methods of the paper in more detail.  Recall that
we have fixed a field $F$, and let $Q_{n-1}\inc \P^{n}$ be the quadric
defined over $F$ by the equation $x_0^2+\cdots+x_n^2=0$.  Write
$DQ_n=\P^n-Q_{n-1}$ (the subscript on a scheme indicates its
dimension).  The scheme $DQ_n$ is in some sense a motivic analogue of
$\RP^n$.  One facet of this analogy is that when
$F=\C$, the space $\RP^n$ is homotopy equivalent to 
the topological space $DQ_n(\C)$ of $\C$-valued points
\cite{Lw}.  
In \cite{SS} and \cite{DI1}
it was observed that if a sums-of-squares formula of type $[r,s,n]$
exists over $F$, then there is an axial map of $F$-varieties
\[ \phi\colon DQ_{r-1}\times DQ_{s-1} \ra DQ_{n-1}.
\]

Here is a list of the main steps in the paper:
\begin{enumerate}[(1)]
\item We prove that if a sums-of-squares formula of type $[r,s,n]$
exists over $\R$ then a diagram of the form
\[ \xymatrix{
\RP^{r-1} \times \RP^{s-1} \ar[d]\ar[r]^{j\times j} &
\RP^\infty \times \RP^\infty \ar[d]^\mu \\
\RP^{n-1}\ar[r]^-j & \RP^\infty 
}
\]
exists in the homotopy category of topological spaces.  Here 
$j$ denotes the obvious inclusions, and $\mu$ is the usual
multiplication map.  Consequently,
if $E^*$ is a complex-oriented cohomology theory then one has 
$(X\tens
1-1\tens X)^e=0$ in $E^*(\RP^{r-1}\times \RP^{s-1})$, where $e=\lfloor
\frac{n+1}{2} \rfloor$.  
We use $X$ to denote 
the complex-orientation class in $E^2(\CP^\infty)$, as well as its
restriction to $E^2(\RP^\infty)$ or to $E^2(\RP^m)$.

\item We next prove that if a sums-of-squares formula of type
$[r,s,n]$ exists over any field $F$ not of characteristic $2$, then
there exists a similar diagram in the homotopy category
$\ho_{\Z/2}(\pro \sSet)$ of pro-simplicial sets and $\Z/2$-cohomological
equivalences.  This uses \'etale homotopy theory and the map $\phi$
mentioned above.  The key step is an analysis of the \'etale homotopy
type of $DQ_n$.

\item 
Suppose $E$ is a complex-oriented cohomology theory that is
$k$-connected for some $k$, and is such that each coefficient group
$E^q$ is a finite $2$-group.  As a formal consequence of the diagram
from (2), one again obtains the identity $(X\tens 1-1\tens X)^e=0$ in
$E^*(\RP^{r-1}\times \RP^{s-1})$.
\item Putting everything together, the existence of a sums-of-squares
formula of type $[r,s,n]$ over an arbitrary field (not of
characteristic $2$) implies that a certain algebraic expression
vanishes in $E^*(\RP^{r-1}\times \RP^{s-1})$, for any cohomology
theory $E$ of the type described in (3).  In particular, one may take
$E^*(X)=BP2^*(X;\Z/2^N)$ for large $N$, in which case one can use the
$BP2$ calculations done by Davis \cite{D2}.  This yields a new
necessary condition for the existence of sums-of-squares formulas over
fields of positive characteristic, as stated in the following theorem.
This result is proved at the end of Section \ref{sec:main}.
\end{enumerate}

\begin{thm}
\label{thm:main}
Suppose a sums-of-squares formula
of type $[2a+1,2b+1,2m-1]$ exists over a field $F$,
where $\max\{a,b\} < m \leq a+b$ and $\chara(F)\neq 2$.  Then the 
vector
\[
[ 
(-1)^a\tbinom{m}{a} , (-1)^{a-1}\tbinom{m}{a-1}, \ldots,  
(-1)^{m-b}\tbinom{m}{m-b}
]
\]
is in the $\Z$-linear span of the relations listed in
Theorem~\ref{thm:DavisBP} below. 
\end{thm}

Unfortunately the relations mentioned in the above theorem (which are
here interpreted as vectors of integers) are too awkward to list in
the introduction.  However, we give the following concrete
application which is discussed in more detail in 
Examples \ref{ex:11-15-17} and \ref{ex:11-15-17-end}.
Consider whether a sums-of-squares formula of type
$[11,15,17]$ exists over a field of characteristic not equal to $2$.
The necessary conditions obtained in \cite{DI1} and \cite{DI2} are
satisfied, which means that those results are inconclusive.  According
to the above theorem, however, if such a formula exists, then the
vector $[-126, 126, -84, 36]$ belongs to the
$\Z$-linear span of 
\[
[2,0,0,0], \quad [0,2,2,2], \quad [0,0,4,4], \quad\text{and}\quad
[0,0,0,2].
\]
A quick calculation shows that a formula of type $[11,15,17]$ 
therefore does not exist.  

\begin{remark}
The paper \cite{BDM} applied the cohomology theory $\tmf$ to the axial
map problem.  Its results can be interpreted as giving further
restrictions on the possible dimensions of sums-of-squares formulas
over $\R$.  According to the authors, those restrictions actually
improve on the $BP2$ restrictions by arbitrarily large amounts as the
dimensions get big.  We have not explored whether these $\tmf$ results
also extend to characteristic $p$ fields (the arguments of the present
paper do not quite apply because $\tmf$ is not complex-oriented); so
this remains an interesting open question.
\end{remark}

\subsection{Organization}
The paper is organized as follows.  Section \ref{sctn:Davis}
is  expository.  It describes the approach of Astey and Davis
to sums-of-squares formulas over $\R$, which is important background
for understanding our arguments.  Section \ref{sec:main}
introduces the basic machinery of our arguments and outlines the
proof of Theorem \ref{thm:main}.  The proofs of two key steps,
Propositions \ref{prop:eta2} and \ref{pr:mainaxialp},
are postponed until later.

The remainder of the paper is dedicated to the task of proving
these two propositions.  Section \ref{sctn:intro2} provides an overview
of the technical issues.
The last three sections fill in the details, as well as provide
more detailed background on the ideas in Section \ref{sec:main}.


\section{Sums-of-squares and complex-oriented theories}
\label{sctn:Davis}

In this section we describe how complex-oriented cohomology theories
can be applied to the sums-of-squares problem over $\R$.  This
technique was developed by Astey \cite{As} and later used by Davis
\cite{D2}.  The specific results we need seem not to be
in the literature, so we briefly develop them here.

\bigskip

Suppose that $E$ is a complex-oriented cohomology theory
\cite[II.2]{Ad} equipped with an orientation class $X$ in
$\tilde{E}^2(\CP^\infty)$.  Recall that complex-oriented cohomology
theories are represented by spectra with a homotopy associative and
homotopy commutative multiplication.

We will abuse notation and also write $X$ for the image of the
orientation class in $\tilde{E}^2(\RP^\infty)$ under the standard
inclusion $\RP^\infty \inc \CP^\infty$.  Even further, we write $X$
also for the image in $\tilde{E}^2(\RP^n)$ under the standard
inclusion $\RP^n \inc \RP^\infty$.  The resulting classes $X\tens 1$
and $1\tens X$ in $E^*(\RP^\infty \times \RP^\infty)$ or $E^*(\RP^m
\times \RP^n)$ will be denoted $X_1$ and $X_2$.

\begin{lemma}
\label{lem:E-power}
Let $E$ be a complex-oriented cohomology theory.
Then $X^e = 0$ in $E^*(\RP^{n-1})$, where $e = \lfloor\frac{n+1}{2}\rfloor$.
\end{lemma}

\begin{proof}
Recall that any complex-oriented cohomology theory $E$ comes
equipped with a (multiplicative) structure 
map of associated spectra $MU \map E$ \cite[II.4.6]{Ad}. 
By naturality, we only need to consider the case $E=MU$.  
Consider the Atiyah-Hirzebruch spectral sequence \cite{AH}
\[
E_2^{p,q} = H^p(\RP^{n-1}; MU^q ) \Rightarrow MU^{p+q}(X).
\]
The $E_2$-term vanishes unless $0 \leq p \leq n-1$ and $q \leq 0$.
It follows that $MU^i(\RP^{n-1})$ is zero for $i>n-1$.

Since the degree of $X^e$ is at least $n$, we conclude that $X^e=0$. 
\end{proof}

The above proof explains the choice of exponent 
$e = \lfloor \frac{n+1}{2} \rfloor$.  It is the smallest integer such
that $2e \geq n$.

\begin{lemma}
\label{lem:Astey}
Let $\mu: \RP^\infty \times \RP^\infty \map \RP^\infty$ be the usual
map that represents tensor product of line bundles,
and let $E$ be a complex-oriented cohomology theory.
Then $\mu^*(X) = u(X_1 - X_2)$, where $u$ is a unit in 
$E^*(RP^\infty\times \RP^\infty)$.
\end{lemma}

\begin{proof}
As in the proof of Lemma \ref{lem:E-power},
it suffices to consider the case $E=MU$, which is 
precisely \cite[Prop.~3.6]{As}.
\end{proof}

\begin{defn}
\label{defn:axial}
A map $\RP^{r-1} \times \RP^{s-1} \map \RP^{n-1}$ is \mdfn{axial} if
there exist basepoints of $\RP^{r-1}$ and $\RP^{s-1}$ such that the
restricted maps $\RP^{r-1} \times * \map \RP^{n-1}$ and $* \times
\RP^{s-1} \map \RP^{n-1}$ take the generator in
$H^1(\RP^{n-1};\Z/2)$ to the generators of $H^1(\RP^{r-1};\Z/2)$
and $H^1(\RP^{s-1};\Z/2)$, respectively.
\end{defn}

\begin{prop}
\label{pr:mainaxial}
Suppose that $f\colon \RP^{r-1} \times \RP^{s-1} \ra \RP^{n-1}$ 
is an axial map,
and let $E$ be a complex-oriented cohomology
theory.  Then $(X_1-X_2)^{e} =0$ in the ring $E^*(\RP^{r-1}\times
\RP^{s-1})$, where $e= \lfloor \frac{n+1}{2} \rfloor$.
\end{prop}

\begin{proof}  
We first claim that there is a homotopy-commutative diagram
\[ \xymatrix{
\RP^{r-1} \times \RP^{s-1} \ar[d]_f\ar[r]^{j\times j} &
\RP^{\infty}\times \RP^\infty \ar[d]^\mu \\
\RP^{n-1}\ar[r]^-j & \RP^\infty.
}
\]
Using that 
$\RP^\infty$ represents $H^1(\blank;\Z/2)$, 
the homotopy-commutativity of this diagram is 
equivalent to the statement that
\[
f^*\colon H^1(\RP^{n-1}; \Z/2) \ra H^1(\RP^{r-1}\times\RP^{s-1}; \Z/2)
\] 
sends the 
generator $x$ to $x \otimes 1 + 1 \otimes x$.  
This is immediate from the
fact that $f$ is an axial map.

Now apply the cohomology theory $E$ to the above diagram.  Together with
Lemma \ref{lem:Astey}, this shows that 
$f^*\colon E^2(\RP^{n-1}) \map E^2(\RP^{r-1} \times \RP^{s-1})$ takes
$X$ to $u(X_1 - X_2)$, where $u$ is a unit in $E^*(\RP^{r-1}\times \RP^{s-1})$.

We also know that $X^e = 0$ in $E^*(\RP^{n-1})$ by Lemma \ref{lem:E-power}.
It follows that $(X_1 - X_2)^e = 0$ in $E^*(\RP^{r-1}\times \RP^{s-1})$.
\end{proof}

Now suppose 
that there is a sums-of-squares formula of type $[r,s,n]$ over $\R$.
This yields a bilinear map $f\colon \R^r \times \R^s \ra \R^n$
given by 
\[
(x_1,\ldots,x_r; y_1,\ldots,y_s) \mapsto (z_1,\ldots,z_n),
\]
using the notation from the first paragraph of the introduction.
This map $f$ has the property that $||f(x,y)||^2=||x||^2 ||y||^2$
for all $x$ and $y$.
In particular, $f(x,y)=0$ only if $x=0$ or $y=0$.  So $f$ induces
a map $\RP^{r-1}\times \RP^{s-1} \ra \RP^{n-1}$, which is readily seen
to be an axial map.  This argument is due to Hopf \cite{H}.
Proposition \ref{pr:mainaxial} therefore gives the following consequence.

\begin{cor}
\label{cor:sos-R}
Suppose that a sums-of-squares formula of type $[r,s,n]$ exists over $\R$.
Then for any complex-oriented cohomology theory $E$, one has
$(X_1-X_2)^e=0$ in the ring $E^*(\RP^{r-1}\times \RP^{s-1})$, where $e=\lfloor
\frac{n+1}{2} \rfloor$.
\end{cor}

For theories $E$ where one can do the relevant
computations, the above corollary leads to necessary conditions for a
formula of type $[r,s,n]$ to exist over $\R$.
The obvious choice is to take $E = MU$.  However,
as Astey \cite{As} points
out, the groups $MU^*(\RP^{r-1}\times \RP^{s-1})$ are very
large---especially compared to 
the subgroup generated by monomials in $X_1$ and $X_2$, which is
what we are really interested in because $(X_1 - X_2)^e$ is a 
linear combination of such monomials.
So it's
preferable to replace $MU$ with something smaller but that still carries
almost as much of the relevant information.  Astey \cite[Sec.~5]{As} used the
Landweber-exact theory
\[ X\mapsto MU^*(X)\tens MU^*/(v_2,v_4,v_5,v_6,\ldots)[v_3^{-1}]\tens
\Z_{(2)},
 \]
for a certain choice of generators $v_1,v_2,\ldots$ of $MU^*$.
Later Davis \cite{D2} superseded these results using the
theory $BP2$.  We will describe Davis's work next.

\subsection{Computations for $BP2$}

Recall that for any prime $p$ there is an associated 
Brown-Peterson cohomology theory $BP$ \cite[II.16]{Ad}.
The coefficients are $BP^*=\Z_{(p)}[v_1,v_2,\ldots]$ with $\deg
v_i=-2(p^i-1)$.  Here $\Z_{(p)}$ is the integers localized at $p$.
Whenever $BP$ is used in the present paper we will have
$p=2$.  The cohomology theory $BP2$ is constructed from
$BP$ by killing the elements $v_3,v_4,\ldots$.
There are at least two ways of constructing $BP2$:
the methods
of Baas-Sullivan \cite{B},
or using that $BP$ is a ring spectrum in a
modern category of spectra \cite{EKMM, HSS}, 
in which case $BP2$ can be constructed as a $BP$-module spectrum 
by iterative cofiber sequences.  
So $BP2^*=\Z_{(2)}[v_1,v_2]$.  The fact that $BP2$ has the structure
of a homotopy associative, homotopy commutative ring spectrum is
classical;  but in fact it is actually an $A_\infty$-spectrum by \cite{An, G}.
Note that $BP2$ is complex-oriented via the composition of the
multiplicative maps $MU \ra MU_{(p)} \ra BP \ra BP2$.

Davis \cite{D2} was able to compute the
groups $BP2^{2m}(\RP^{2a}\times \RP^{2b})$ completely for $m$ in a
range.  The following result is \cite[Thm.~1.4]{D2}:

\begin{thm}[Davis] 
\label{thm:DavisBP}
Assume that $\max\{a,b\}< m \leq a+b$; let $d=a+b-m$.  Then
$BP2^{2m}(\RP^{2a}\times \RP^{2b})$ is the abelian group
generated by the $d+1$ elements
\[ w_0 = X_1^a X_2^{m-a}, w_1 = X_1^{a-1}X_2^{m-a+1}, \ldots, 
w_i = X_1^{a-i}X_2^{m-a+i}, \ldots, w_d = X_1^{m-b}X_2^b
\]  
subject to the $d+1$ relations:
\begin{enumerate}[(a)]
\item 
$2w_0 =0$.

\vspace{0.1in}

\item 
$\displaystyle 2^{\lfloor \frac{i+2}{2} \rfloor} 
\left[ \tbinom{c}{c} w_i + \tbinom{c+1}{c} w_{i+1} + 
\tbinom{c+2}{c} w_{i+2} + \cdots + \tbinom{c+d-i}{c} w_d  \right]=0$ for
$c=\lfloor \frac{i-1}{2} \rfloor$ and
 $1\leq i < d+1-\bigl\lfloor{\tfrac{d+1}{3}}\bigr \rfloor$.

\vspace{0.1in}

\item $2^{d+1-i}w_i =0$ 
if $d+1-\bigl\lfloor{\tfrac{d+1}{3}}\bigr \rfloor \leq i \leq d$.
\end{enumerate}
\end{thm}

To better understand the above relations, it is helpful to 
explore some examples.  

\begin{ex}
Consider $BP2^{22}(\RP^{20} \times \RP^{20})$.
Then $a=b=10$, $m=11$, and $d=9$.  The group is the
quotient of $\Z^{10}$ by the rowspace of the matrix
\[
\begin{bmatrix}
2 & 0 & 0 & 0 & 0 & 0 & 0 & 0 & 0 & 0 \\
0 & 2 & 2 & 2 & 2 & 2 & 2 & 2 & 2 & 2 \\
0 & 0 & 4 & 4 & 4 & 4 & 4 & 4 & 4 & 4 \\
0 & 0 & 0 & 4 & 8 & 12 & 16 & 20 & 24 & 28 \\
0 & 0 & 0 & 0 & 8 & 16 & 24 & 32 & 40 & 48 \\
0 & 0 & 0 & 0 & 0 & 8 & 24 & 48 & 80 & 120 \\
0 & 0 & 0 & 0 & 0 & 0 & 16 & 48 & 96 & 160 \\
0 & 0 & 0 & 0 & 0 & 0 & 0 & 8 & 0 & 0 \\
0 & 0 & 0 & 0 & 0 & 0 & 0 & 0 & 4 & 0 \\
0 & 0 & 0 & 0 & 0 & 0 & 0 & 0 & 0 & 2
\end{bmatrix}.
\]
The general pattern can be described easily.
Corresponding to relation (a),
the first row always consists of a $2$
followed by zeros.  
Corresponding to the relations in (c),
the bottom-right portion of the matrix  is always
diagonal, with consecutive powers of $2$ ascending on the diagonal
from the bottom-right corner; there are  precisely $\lfloor
\frac{d+1}{3}\rfloor$ of these diagonal entries.  
Corresponding to the relations in (b), the rest
of the matrix can be filled in as follows.  Consider the
sequences
\begin{align*}
 & \ 1,1,1,\ldots \\
 &\tbinom{1}{1},\tbinom{2}{1},\tbinom{3}{1},\ldots \\
 &\tbinom{2}{2},\tbinom{3}{2},\tbinom{4}{2},\ldots \\
 & \tbinom{3}{3},\tbinom{4}{3},\tbinom{5}{3},\ldots 
\end{align*}
and so forth.
Each of these sequences appears in two rows of the matrix (given
enough room), multiplied by consecutive powers of $2$, and shifted so
that the first entry appears on the diagonal.
\end{ex}

\begin{ex}
\label{ex:11-15-17}
Consider the group $BP2^{18}(\RP^{10}\times
 \RP^{14})$.  Here $a= 5$, $b=7$, and $m=9$, so $d=3$ and 
$\lfloor \frac{d+1}{3} \rfloor =1$.
The group is the quotient of $\Z^4$ by the rowspace of
\[ \begin{bmatrix}
2 & 0 & 0 & 0 \\
0 & 2 & 2 & 2 \\
0 & 0 & 4 & 4 \\
0 & 0 & 0 & 2
\end{bmatrix}.
\]
\end{ex}

\begin{cor}
\label{cor:davis-char=0}
Suppose that there exists a sums-of-squares formula over $\R$ 
of type $[2a+1,2b+1,2m-1]$,
where $\max\{a,b\} < m \leq a+b$.  Then the vector
\[ [(-1)^a\tbinom{m}{a}, (-1)^{a-1}\tbinom{m}{a-1},\ldots,
(-1)^{m-b}\tbinom{m}{m-b}]
\]
is in the $\Z$-linear span of the relations listed in
Theorem~\ref{thm:DavisBP} (where we interpret the relations as vectors
of coefficients). 
\end{cor}

\begin{proof}
From Corollary~\ref{cor:sos-R}, $(X_1 - X_2)^m = 0$ in 
$BP2^{2m}(\RP^{2a} \times \RP^{2b})$ because
$m = \left\lfloor \frac{2m-1+1}{2} \right\rfloor$.  Expand this expression
to obtain
\[
\sum_{k=0}^{m} (-1)^k \tbinom{m}{k} X_1^k X_2^{m-k} = 0.
\]
By Lemma~\ref{lem:E-power},
$X_1^k X_2^{m-k} = 0$ unless $m-b \leq k \leq a$.  So we can simplify
the above equation to
\[
\sum_{k=m-b}^a (-1)^k \tbinom{m}{k} w_{a-k} = 0;
\]
here we are also using that $w_{a-k} = X_1^k X_2^{m-k}$
if $m-b \leq k \leq a$.
\end{proof}

\begin{ex}
\label{ex:11-15-17-end}
Consider the question of whether a sums-of-squares formula of type
$[11,15,17]$ exists over $\R$.  Apply Corollary \ref{cor:davis-char=0} with
$a=5$, $b=7$, and $m=9$.  
According to Example \ref{ex:11-15-17},
we need to check whether the
vector
$
[-\tbinom{9}{5}, \tbinom{9}{4},-\tbinom{9}{3}, \tbinom{9}{2} ]
$
is in the $\Z$-linear span of the vectors
\[
[2,0,0,0], \quad [0,2,2,2], \quad [0,0,4,4], \quad\text{and}\quad
[0,0,0,2].
\]
One readily checks that it's {\it not\/} in the span, 
so formulas of type
$[11,15,17]$ do not exist over $\R$.  
Note that this trivially implies that a sums-of-squares formula 
of type $[11,15,k]$
does not exist for any $k\leq 17$.

By comparison, the Hopf-Stiefel condition \cite{H, St}
only shows that
formulas of type $[11, 15, 15]$ do not exist, while the Atiyah
condition \cite{A} only shows that formulas of type $[11, 15, 16]$ (or
equivalently, $[15,11,16]$) do not exist.
\end{ex}


\section{The main results}
\label{sec:main}

This section gives the main results of the paper.  We begin by
summarizing the technical machinery needed.

\subsection{Pro-spaces}
\label{subsctn:pro}

Let \mdfn{$\pro\sSet$} denote the category of pro-simplicial sets.  Recall that
an object of this category consists of a small, cofiltered category $I$
together with a functor $X\colon I \ra \sSet$.  If $Y\colon J \ra
\sSet$ is another such diagram, then the morphisms from $X$ to $Y$ are
defined to be
\[ \Hom_{\pros}(X,Y)=\lim_j \colim_i \Hom_{\sSet}(X_i,Y_j).
\]
Objects of $\sSet$ will be called `spaces', and so objects of $\pro \sSet$
are `pro-spaces'.
Traditional references for pro-categories include
\cite{AM} and \cite[Sec.~I.8]{SGA4}, but for our purposes \cite{I1}
is more relevant.

If $X$ is a pro-space and $M$ is an abelian group, then the
\mdfn{singular cohomology $H^*(X;M)$} of $X$ with coefficients in $M$
is defined to be $\colim_s H^*(X_s;M)$ \cite[2.2]{AM}.  Since $X$ is a
cofiltered system and cohomology is contravariant, the colimit in the
definition of $H^*(X;M)$ is a filtered colimit.

A \mdfn{constant pro-space} is one indexed by the trivial category
with one object and an identity map.
There is a ``constant'' functor \mdfn{$c$}$\colon\sSet \ra \pro \sSet$
that takes a simplicial set $X$ to the constant pro-space with value $X$.
Note that for any space $X$ and any abelian group $M$, the cohomology rings
$H^*(X;M)$ and $H^*(cX;M)$ are naturally isomorphic.

We will frequently abuse notation in the following way.
If $X$ is a topological space, we write $cX$
for the constant functor applied to
the singular simplicial set of $X$.

There are several useful model category structures on pro-spaces.  The
one that is relevant for us is the $\Z/2$-cohomological model structure
\cite{I3}.  
A map of pro-spaces $X\ra Y$ is called a
\mdfn{$\Z/2$-cohomological equivalence} if the induced map $H^*(Y;\Z/2)\ra
H^*(X;\Z/2)$ is an isomorphism.  
The presence of the model structure establishes the existence of
the associated homotopy category
$\ho_{\Z/2}(\pro \sSet)$ obtained by formally inverting the class
of $\Z/2$-cohomological equivalences.

We will need one additional property of the $\Z/2$-cohomological
homotopy theory of pro-spaces.  The proof of the following result
will be postponed until Section \ref{sctn:pro}.

\begin{prop}
\label{prop:eta2}
For any spectrum $E$, there is a naturally defined functor
\[ (\Ecomp_2)^*\colon\pro \sSet \ra (\text{graded abelian groups})  
\]
called the \mdfn{$2$-completion of $E$}, satisfying
the following properties:
\begin{enumerate}[(i)]
\item $(\Ecomp_2)^*$ takes $\Z/2$-cohomological equivalences to
  isomorphisms;
\item There is a natural transformation $\eta\colon E^*(X) \ra
(\Ecomp_2)^*(cX)$ for all spaces $X$;  
\item If $E$ is $k$-connected for some
$k$ and each coefficient group $E^q$ is a finite $2$-group,
then $\eta$ is a natural isomorphism. 
\end{enumerate}
\end{prop}

\subsection{The main argument}
\label{subsctn:main}
The following proposition is the central technical ingredient in the proof 
of Theorem \ref{thm:main}, our main result 
about sums-of-squares formulas.
Its proof will be postponed until the end of Section \ref{sctn:et}.

\begin{prop}
\label{pr:mainaxialp}
Let $F$ be a field of characteristic not equal to $2$, and
suppose a sums-of-squares formula of type $[r,s,n]$ exists over $F$.
Then there exists a commutative diagram in $\ho_{\Z/2}(\pro \sSet)$ of the form
\[\xymatrixcolsep{3.6pc}
\xymatrix{ 
c(\RP^{r-1}\times \RP^{s-1}) \ar[d]
   \ar[r]^{c(j\times j)} 
& 
c(\RP^\infty \times \RP^\infty) \ar[d]^{c(\mu)} \\
c\RP^{n-1} \ar[r]^{c(j)} & c\RP^\infty.  }
\]
Here each $j$ is the standard inclusion, and
$\mu\colon \RP^\infty \times\RP^\infty\ra \RP^\infty$ is the map
classifying tensor product of line bundles.
\end{prop}

\begin{remark}
The reader should compare the above diagram to the one from the proof
of Proposition~\ref{pr:mainaxial}.  
We do not claim that the left vertical map is induced by a map
of ordinary spaces.  It does exist in $\ho_{\Z/2}(\pro \sSet)$,
but it is not necessarily in the image of the constant functor $c$
because it is constructed using a homotopy-inverse in 
$\ho_{\Z/2}(\pro \sSet)$; see the second diagram in the proof
of Proposition \ref{pr:mainaxialp} given in Section \ref{sctn:et}.
\end{remark}

We will give the proof of the above proposition in subsequent
sections; this is where \'etale homotopy theory enters the picture.
For the moment, however, we will just accept the proposition and
consider its consequences.

\begin{prop}
\label{prop:axialp}
Let $E$ be a complex-oriented cohomology theory that is $k$-connected
for some $k$.  Suppose that
a sums-of-squares formula of type $[r,s,n]$ exists over a field $F$ of
characteristic not equal to $2$.  Then $(X_1-X_2)^e=0$ in
$E^*(\RP^{r-1}\times \RP^{s-1})$, where $e=\lfloor
\frac{n+1}{2}\rfloor$, provided either of the following two conditions
is satisfied:
\begin{enumerate}[(a)]
\item
Each $E^q$ is a finite $2$-group;
\item  The
kernel and cokernel of $E^{q}\llra{2} E^{q}$ are
finite-dimensional $\Z/2$-vector spaces for each $q$, 
and both of $r$ and $s$ are odd.
\end{enumerate}
\end{prop}

Several small variations of the above result are possible.  For
instance, if one assumes that $E^*$ is concentrated in even degrees
then one can prove (b) assuming only that one of $r$ and $s$ is odd.
\begin{proof}
In case (a) we begin by applying the functor $(\Ecomp_2)^*$ to the
diagram from Proposition~\ref{pr:mainaxialp} to obtain a diagram of
graded abelian groups.  We also have the diagram
\[\xymatrixcolsep{3pc}
\xymatrix{ 
\RP^{r-1}\times \RP^{s-1} \ar[r]^{j \times j} & 
  \RP^\infty \times \RP^\infty \ar[d]^\mu \\
\RP^{n-1} \ar[r]_j & \RP^\infty     }
\]
of spaces, to which we apply the functor $E^*$.  The natural
transformation $\eta$ from Proposition \ref{prop:eta2} gives us a commutative
diagram
\[
\xymatrix@C-10ex{
& (\Ecomp_2)^*c(\RP^{r-1} \times \RP^{s-1}) & & 
            (\Ecomp_2)^*c(\RP^\infty \times \RP^\infty) \ar[ll] \\
E^*(\RP^{r-1} \times \RP^{s-1}) \ar[ur] & & 
            E^*(\RP^\infty \times \RP^\infty) \ar[ll]\ar[ur] \\
& (\Ecomp_2)^* (c\RP^{n-1}) \ar[uu]|(.5)\hole & & 
            (\Ecomp_2)^* (c\RP^\infty) \ar[uu]\ar[ll]|(.5)\hole \\
E^* (\RP^{n-1}) \ar[ur] & & E^* (\RP^\infty) \ar[uu]\ar[ll]\ar[ur]   }
\]
in which all four diagonal maps are isomorphisms.
Lemma \ref{lem:E-power} shows that 
the element $X^e$ of $E^*(\RP^\infty)$ maps to zero in
$(\Ecomp_2)^* c(\RP^{r-1} \times \RP^{s-1})$.
Because $\eta$ is an isomorphism, $X^e$ also maps to
zero in $E^*(\RP^{r-1} \times \RP^{s-1})$.  
On the other hand, 
Lemma \ref{lem:Astey}
says that $X^e$ maps to a unit multiple of
$(X_1-X_2)^e$ in 
$E^*(\RP^{r-1}\times \RP^{s-1})$.  
We have computed the image of $X^e$ in two different ways, and
we conclude that $(X_1 - X_2)^e = 0$.

Case (b) requires more work.  
Recall that the reduced $E$-cohomology $\tilde{E}^*(Z)$ of a 
pointed space $Z$ is defined to be the kernel of the map
$E^*(Z) \map E^*$ induced by $* \map Z$.  By the definition of $X_1$
and $X_2$ in terms of the orientation class of the complex-oriented
theory $E$, we know that $(X_1 - X_2)^e$ belongs to
$\tilde{E}^*(\RP^{r-1} \times \RP^{s-1})$.  Therefore, we are allowed
to compute in reduced $E$-cohomology instead of unreduced $E$-cohomology.

We examine the
reduced Atiyah-Hirzebruch spectral sequence \cite{AH}
\[
E^{p,q}_2 = \tilde{H}^p ( \RP^{r-1} \times \RP^{s-1}; E^q ) 
  \Rightarrow \tilde{E}^{p+q} ( \RP^{r-1} \times \RP^{s-1} ).
\]
If $r$ and $s$ are odd, then the singular homology groups
$\tilde{H}_*(\RP^{r-1}\times \RP^{s-1};\Z)$ are all
finite-dimensional $\Z/2$-vector spaces and vanish above dimension
$r+s-2$.  
Our assumptions on $E^*$ guarantee that $\Ext_{\Z}(V, E^q)$
and $\Hom_{\Z}(V, E^q)$ are finite-dimensional $\Z/2$-vector spaces
for any finite-dimensional $\Z/2$-vector space $V$.
The universal coefficient theorem 
shows that each group $E_2^{p,q}$ in our
spectral sequence is also a finite-dimensional $\Z/2$-vector space,
and that these groups vanish when $p>r+s-1$.  
It follows immediately
that $\tilde{E}^k(\RP^{r-1} \times \RP^{s-1})$ is a finite $2$-group
for all $k$.  Thus, there exists a sufficiently large $N$ such that
$2^N$ kills $\tilde{E}^{2e}(\RP^{r-1}\times\RP^{s-1})$.

Consider the cohomology theory $E^*(\blank;\Z/2^N)$.  As long as
$N\geq 3$, the spectrum representing this theory
inherits a homotopy associative and commutative
multiplication from $E$ \cite{O}.  So the theory
is complex-oriented.  
Moreover, the theory satisfies the hypotheses of case (a).
So we have that $(X_1-X_2)^e=0$ in $\tilde{E}^{2e}(\RP^{r-1}\times
\RP^{s-1};\Z/2^N)$.  But consider the following three terms of a long
exact sequence
\[  \tilde{E}^{2e}(\RP^{r-1}\times \RP^{s-1}) \llra{2^N}
\tilde{E}^{2e}(\RP^{r-1}\times \RP^{s-1}) \ra
\tilde{E}^{2e}(\RP^{r-1}\times \RP^{s-1};\Z/2^N).
\]
The first map is zero by our choice of $N$, so
the second map is a monomorphism.
Therefore, $(X_1-X_2)^e=0$ in $\tilde{E}^{2e}(\RP^{r-1}\times \RP^{s-1})$.
\end{proof}

Proposition \ref{prop:axialp}
lets one adapt many known results in
characteristic zero to characteristic $p>2$, simply by varying the
choice of $E$.  The cohomology theory
$MU/(v_2,v_4,v_5,v_6,\ldots)[v_3^{-1}]\tens \Z_{(2)}$ used by Astey
\cite {As} can serve as one such choice, 
but its results are eclipsed
by those of $BP2$.

\begin{thm}
\label{thm:davis-char=p}
Let $F$ be a field of characteristic not equal to $2$.  Suppose that a
sums-of-squares formula of type $[2a+1,2b+1,2m-1]$ exists over $F$,
where $\max\{a,b\}< m \leq a+b$.  Then the vector
\[
\bigl [\ (-1)^a\tbinom{m}{a}, \ (-1)^{a-1}\tbinom{m}{a-1},\  \ldots,\  
(-1)^{m-b}\tbinom{m}{m-b}\ \bigr ]
\]
is in the $\Z$-linear span of the relations listed in
Theorem~\ref{thm:DavisBP}. 
\end{thm}

\begin{proof}
Apply Proposition~\ref{prop:axialp} with $E=BP2$.  Note that $BP2^*$ is
$\Z_{(2)}[v_1, v_2]$, where 
$v_1$ and $v_2$ are classes in dimensions
$-2$ and $-6$ respectively, and $\Z_{(2)}$ is the integers localized
at 2.  
It is immediate that $BP2$ satisfies
condition (b) of Proposition \ref{prop:axialp}.  
Then use Davis's computations
from Theorem~\ref{thm:DavisBP}, just as in the proof of Corollary
\ref{cor:davis-char=0}.
\end{proof}

\begin{remark}
In this section we have used the $\Z/2$-cohomological model structure
on pro-simplicial sets.  We strongly suspect that the homotopy theory
of \cite{Mo} is also a suitable context for the above results.
\end{remark}


\section{Introduction to the rest of the paper}
\label{sctn:intro2}

The remaining sections of the paper are concerned with the proofs of
Propositions \ref{prop:eta2} and \ref{pr:mainaxialp}.  These are
the two uncompleted parts of the main argument given in 
Section \ref{sec:main}.

\medskip

We will be primarily 
concerned with the following specific schemes.

\begin{defn}
\label{defn:DQ}
If $R$ is a ring, let
\mdfn{$Q_{n-1}$} be the quadric hypersurface in $\P^n_R$ defined by
the equation $x_0^2 + \cdots + x_n^2 = 0$.
Let
\mdfn{$DQ_n$} be the open complement $\P^n_R - Q_{n-1}$.
\end{defn}

Usually the ring $R$ will be clear from context, but in cases where it
isn't we will write \mdfn{$Q_{n-1/R}$} and \mdfn{$DQ_{n/R}$}.

\begin{remark}
\label{rem:complex-case}
For any scheme defined over $\C$, let $X(\C)$ denote the set of
$\C$-valued points of $X$ equipped with the complex-analytic topology.
From \cite[Lem.~6.3]{Lw}, we know that $DQ_n(\C)$ is homotopy
equivalent to $\RP^n$.  So $(DQ_m \times DQ_n)(\C)$ is homotopy
equivalent to $\RP^m \times \RP^n$.
\end{remark}

Here is a brief
outline of the remaining steps in the paper:
\begin{enumerate}[(1)]
\item 
In Section \ref{sctn:pro} we prove Proposition \ref{prop:eta2}.  We study
generalized cohomology theories and the $\Z/2$-cohomological model
stucture for pro-spaces in more detail than in Section
\ref{subsctn:pro}.
\item 
Section \ref{sctn:et-cohlgy} concerns the \'etale cohomology of
the schemes $DQ_n$ and $DQ_m \times DQ_n$
(see Examples \ref{ex:H-DQ} and \ref{ex:H-DQxDQ}).
We will show in Theorem \ref{thm:P-linear-cohomology}
that the \'etale $\Z/2$-cohomology rings of these schemes
do not change under certain kinds of extensions of scalars;
in fact, they are 
isomorphic to the singular $\Z/2$-cohomology rings of $\RP^n$
and $\RP^m \times \RP^n$.
In order to do this, we first have to extend some standard results about
\'etale cohomology from the class of smooth, projective schemes 
to the larger class
of ``$\SP$-linear'' schemes.  These are the schemes that can be
built from smooth, projective schemes by iteratively taking
complements of closed inclusions.  The ideas in this section are
entirely algebro-geometric, not homotopical.

\item 
In Section~\ref{sctn:et} we will study the \'etale realizations
of $\DQ_n$ and $\DQ_m \times \DQ_n$.  We will show that over a 
separably closed field these \'etale realizations are
$\Z/2$-cohomologically equivalent to the constant pro-spaces
$c\RP^n$ and $c(\RP^m \times \RP^n)$
(see Corollary \ref{cor:2-wk-eq}).
This allows us to complete the proof of Proposition~\ref{pr:mainaxialp}.
\end{enumerate}


\section{Generalized cohomology of pro-spaces}
\label{sctn:pro}

In this section we expand our discussion of pro-spaces
from Section \ref{subsctn:pro}.  
We define the 2-complete $E$-cohomology groups $(\Ecomp_2)^* X$
when $X$ is a pro-space and $E$ is an ordinary spectrum.  
We then prove Proposition \ref{prop:eta2}.
We assume that the reader is familiar with the language of 
model structures \cite{Hi, Ho, Q}.

\subsection{$\Z/2$-cohomological homotopy theory of pro-spaces}

Recall that there is a $\Z/2$-cohomological model structure
on the category of pro-spaces.  The $\Z/2$-cohomological equivalences
are the pro-maps that induce isomorphisms on singular $\Z/2$-cohomology.
We write \mdfn{$[\blank, \blank]^{\pros}_{\Z/2}$} for sets of homotopy classes
in the $\Z/2$-cohomological homotopy category $\ho_{\Z/2}(\pro \sSet)$.

The cofibrations of pro-spaces are easy to describe.  They are the
pro-maps that are isomorphic to levelwise cofibrations.  In particular,
every pro-space is cofibrant because every simplicial set is cofibrant.

The $\Z/2$-cohomologically fibrant objects are 
described by \cite[Thm. 3.3]{I3}, but the complete description is
complicated and more than we actually need here.  We need to know 
that the constant pro-space $cK(\Z/2,n)$ is 
$\Z/2$-cohomologically fibrant if
$K(\Z/2,n)$ is a fibrant Eilenberg-Mac Lane space \cite[Thm.~3.3]{I3}.
In particular, $c\RP^\infty$ is $\Z/2$-cohomologically fibrant.  It follows
formally that the product 
$c\RP^\infty \times c\RP^\infty = c(\RP^\infty \times \RP^\infty)$
is also $\Z/2$-cohomologically fibrant.

\begin{lemma}
\label{lem:rep-H1}
For any pro-space $X$, $[X, c\RP^\infty]^{\pros}_{\Z/2}$ 
is naturally isomorphic
to $H^1(X;\Z/2)$, and $[X, c(\RP^\infty \times \RP^\infty)]^{\pros}_{\Z/2}$
is naturally isomorphic to $H^1(X;\Z/2) \times H^1(X;\Z/2)$.
\end{lemma}

\begin{proof}
We just observed
that both $c\RP^\infty$ and $c(\RP^\infty \times \RP^\infty)$
are $\Z/2$-co\-hom\-olog\-ically fibrant.
Since $X$ is automatically cofibrant,
we know that $[X,c\RP^\infty]^{\pros}_{\Z/2}$
is equal to 
\[
\pi_0 \Map_{\pros}(X, c\RP^\infty) = 
\pi_0 \colim_s \Map(X_s, c\RP^\infty).
\]
Using that $\pi_0$ commutes with filtered colimits, this last
expression can be identified with 
$\colim_s H^1(X_s;\Z/2)=H^1(X;\Z/2)$.

This establishes the first statement. The proof of the second is similar.
\end{proof}

\begin{remark}
\label{rem:H1}
If $X$ is an ordinary space, then the lemma implies that the
maps
$[X, \RP^\infty] \map [cX, c\RP^\infty]^{\pros}_{\Z/2}$ and
$[X, \RP^\infty \times \RP^\infty] \map 
  [cX, c(\RP^\infty \times \RP^\infty)]^{\pros}_{\Z/2}$
induced by the functor $c$ are both isomorphisms.
\end{remark}

\subsection{$\Z/2$-cohomological model structure for pro-spectra}

In addition to the $\Z/2$-cohomological model structure on
$\pro \sSet$, there is an analogous $\Z/2$-cohomological model
structure on the category of pro-spectra.  
As for pro-spaces, the weak equivalences
are the pro-maps that induce isomorphisms in singular cohomology with
coefficients in $\Z/2$.
The existence of the $\Z/2$-cohomological model structure on pro-spectra
can be established with the general result \cite[Thm.~2.2]{I3}, where $K$
is a set of fibrant models for the
Eilenberg-Mac Lane spectra $\Sigma^{-q} H\Z/2$.
We write \mdfn{$[\blank, \blank]^{\pros}_{\Z/2}$} for sets of
homotopy classes in the associated $\Z/2$-cohomological 
homotopy category of pro-spectra.

As for pro-spaces, the cofibrations are easy to describe.  They are the
pro-maps that are isomorphic to levelwise cofibrations.  In particular,
a pro-spectrum $X$ is cofibrant if each $X_s$ is a cofibrant spectrum.

\begin{defn}
\label{defn:finite-spectra}
A spectrum $X$ is \mdfn{$\Z/2$-finite} if the stable homotopy groups
$\pi_i X$ are all finite 2-groups and are zero for all but finitely
many values of $i$.
\end{defn}

The $\Z/2$-cohomologically fibrant pro-spectra are
described by a stable version of \cite[Thm.~3.3]{I3}, 
but again the complete description is complicated
and more than we actually need here.  
For our purposes, it suffices to observe that a pro-spectrum $X$ is
$\Z/2$-cohomologically fibrant if certain structure maps associated to
the diagram $X$ are fibrations and if each $X_s$ is $\Z/2$-finite.  The
idea of the proof is to use the following lemma to simplify
the notion of $K$-nilpotence from \cite[Defn.~3.1]{I3}.

\begin{lemma}
\label{lem:finite}
Let $\cC$ be the smallest class of spectra such that:
\begin{enumerate}
\item
$\cC$ contains the trivial spectrum $*$.
\item
If $X \map Y$ is a weak equivalence, then $X$ belongs to 
$\cC$ if and only if $Y$ belongs to $\cC$.
\item
If $\Sigma^{-q} H\Z/2 \map X \map Y$ is a 
homotopy fiber sequence and $Y$ belongs to $\cC$, then so
does $X$.
\end{enumerate}
The class $\cC$ is equal to the class of $\Z/2$-finite spectra.
\end{lemma}

\begin{proof}
We first show that $\cC$ is contained in the class of $\Z/2$-finite spectra.
To do this,
note first that the trivial spectrum $*$ is $\Z/2$-finite.
Second, if $X \map Y$ is a weak equivalence, then $X$ is $\Z/2$-finite
if and only if $Y$ is $\Z/2$-finite.  Third, the long exact sequence
of homotopy groups shows that if $\Sigma^{-q} H\Z/2 \map X \map Y$
is a homotopy fiber sequence and $Y$ is $\Z/2$-finite, then 
so is $X$.  This finishes one direction.

For the other direction, let $X$ be a $\Z/2$-finite spectrum.  We want
to show that $X$ belongs to $\cC$.  The proof is by induction on the
order of the finite group $\oplus \pi_i X$.
The base case is that contractible spectra belong to $\cC$.

Let $n$ be the largest number such that $\pi_n X$ is non-zero.  Choose
a short exact sequence
\[
0 \map \Z/2 \map \pi_n X \map A \map 0.
\]
This is always possible because $\pi_n X$ is a finite abelian 2-group.
One can construct a homotopy fiber sequence
\[
\Sigma^n H\Z/2 \map X \map X',
\]
where $\pi_i X'$ equals $\pi_i X$ for $i \neq n$ 
and $\pi_n X'$ equals $A$.
Note that $X'$ is $\Z/2$-finite, and
the order of $\oplus \pi_i X'$ is strictly
smaller than the order of $\oplus \pi_i X$.  By the induction assumption,
$X'$ belongs to $\cC$.
Using property (3) of the class $\cC$,
the fiber sequence now shows that $X$ also belongs to $\cC$.
\end{proof}

\begin{lemma}
\label{lem:cQ}
The functor $c$ is a left Quillen functor from the ordinary stable model
structure on spectra to the $\Z/2$-cohomological model structure
on pro-spectra.
\end{lemma}

\begin{proof}
By definition of pro-categories, $c$ is left adjoint to the functor $\lim$
that takes a pro-object $X$ to the cofiltered limit $\lim_s X_s$ of the
diagram $X$.  Note that $c$ takes cofibrations to levelwise cofibrations,
and $c$ takes weak equivalences to $\Z/2$-cohomological weak equivalences
because ordinary stable weak equivalences 
induce isomorphisms in $\Z/2$-cohomology.
\end{proof}

\subsection{2-complete cohomology theories}

For any pro-spectrum $X$ and any $q\in \Z$, the pro-spectrum
$\Sigma^{q} X$ is constructed by applying the functor $\Sigma^{q}$
to each spectrum $X_s$ in the diagram $X$.

To simplify notation, we will always assume that $E$ is a fibrant
spectrum (in any of the standard simplicial model category structures
for spectra).  Also, we assume that Postnikov sections $P_n E$ have been
defined functorially in such a way that the natural maps $P_n E \map
P_{n-1} E$ in the Postnikov tower of $E$ are fibrations between
fibrant spectra.

The $\Z/2$-cohomological model structure for pro-spectra allows us
to construct 2-complete cohomology theories.

\begin{defn}
\label{defn:2-comp}
Let $E$ be a spectrum, and let $X$ be a pro-spectrum.  The groups
\mdfn{$(\Ecomp_2)^q(X)$} are defined to be $[\Sigma^{-q}X, cE]_{\Z/2}^{\pros}$.
For any spectrum $X$, 
the natural map \mdfn{$\eta$}$: E^q(X) \map (\Ecomp_2)^q(cX)$ is the
map $[\Sigma^{-q} X, E] \map [c\Sigma^{-q} X, cE]_{\Z/2}^{\pros}$
induced by the functor $c$.
\end{defn}

Lemma \ref{lem:cQ} explains why $\eta$ is well-defined.
See \cite{Qk} for a different approach to the same basic construction.

If $X$ is a
pro-space, we write $(\Ecomp_2)^*(X)$ for $(\Ecomp_2)^*(\Sigma^\infty X_+)$.  
Here
$\Sigma^{\infty} X_+$ is constructed by applying the functor $Z\mapsto
\Sigma^\infty(Z_+)$ to each space $X_s$ in the cofiltered diagram $X$.

In order to compute $(\Ecomp_2)^*(X)$, one needs to take a 
$\Z/2$-cohomologically fibrant replacement for $cE$.
Next we consider situations in which this fibrant replacement
and the map $\eta$ are computable.

\begin{lemma}
\label{lem:fibrant}
Suppose that $E$ is a spectrum 
such that $E$ is $k$-connected for some $k$ and each coefficient group
$E^q$ is a finite $2$-group.  
The Postnikov tower $PE$ of $E$, viewed as a pro-spectrum,
is a $\Z/2$-cohomological fibrant replacement for $cE$.
\end{lemma}

\begin{proof}
To begin, note that the map $cE \map PE$ induces an isomorphism
in $\Z/2$-cohomology because $E \map P_n E$ is an isomorphism
in cohomology in degrees less than $n$.  
It remains only to
observe that each $P_n E$ is $\Z/2$-finite.
\end{proof}

Finally, we come to the main point of this section, which is to prove
Proposition~\ref{prop:eta2}.  

\begin{proof}[Proof of Proposition~\ref{prop:eta2}]
For part (i), note that $(\Ecomp_2)^*$ takes $\Z/2$-cohomological equivalences
to isomorphisms because $(\Ecomp_2)^*$ is defined in terms
of the $\Z/2$-cohomological homotopy category.

For part (ii),
we constructed the natural transformation $\eta$ in
Definition \ref{defn:2-comp}.

Finally, for part (iii), 
let $E$ be a spectrum that is $k$-connected for some
$k$.  Also suppose that each coefficient group $E^q$ is a finite $2$-group.
We may assume that $\Sigma^{-q} X$
is a cofibrant spectrum.
Lemma \ref{lem:fibrant} says that
the Postnikov tower $PE$ is a $\Z/2$-cohomological fibrant replacement 
for $cE$.  
Then $\Map_{\pros}(c\Sigma^{-q}X, PE)$ is equal
to $\Map(\Sigma^{-q}X, \lim_n P_n E)$ by adjointness, which
is weakly equivalent to $\Map(\Sigma^{-q}X, E)$
because $E$ is weakly equivalent to the limit of its own Postnikov tower.
Now apply $\pi_0$ to these mapping spaces to show that 
$\eta$ is an isomorphism.
\end{proof}


\section{\'Etale cohomology and $\SP$-linear schemes}
\label{sctn:et-cohlgy}

Recall the schemes $Q_{n-1}$ and $DQ_n$ defined in Section \ref{sctn:intro2}.
We wish to compute the \'etale $\Z/2$-cohomology rings
of the schemes $DQ_n$ and $DQ_m \times DQ_n$ over certain base rings,
and we wish to understand the behavior of these rings under certain
kinds of extensions of scalars.
These schemes are smooth but not proper.  The technical details
of this section get around the absence of properness.

In this section, we are working entirely with ideas from algebraic geometry.
Homotopy theory will enter the picture again in Section \ref{sctn:et} when 
we study the \'etale realizations of $DQ_n$ and $DQ_m \times DQ_n$.

We begin with an issue of notation.  If $X$ is a scheme defined over
$\Spec R$ and $R \map S$ is any ring map, then \mdfn{$X_S$} is the scheme
$\Spec S \times_{\Spec R} X$.  This notation will appear throughout
this and the next section.

\begin{defn}
\label{defn:P-linear}
Let $S$ be a Noetherian scheme.
Let \mdfn{$\SP$} be the class of schemes that are smooth and proper
over $S$.  The class of \mdfn{$\SP$-linear schemes over $S$} 
is the smallest class
of smooth schemes over $S$ containing $\SP$ such that if $Z \map X$
is a closed inclusion of smooth schemes with open complement $U$ and
any two of $Z$, $U$, and $X$ are $\SP$-linear, then so is the third.
\end{defn}

Recall that the linear schemes \cite{J} are the smallest class of
schemes that contain the affine spaces $\A_S^n$ and are also closed
under the same two-out-of-three property.  This motivates the
terminology of Definition \ref{defn:P-linear}.
In fact, $\A_S^n$ is $\SP$-linear because it is the
complement $\P^n_S-\P^{n-1}_S$.  So every linear scheme is
$\SP$-linear.

\begin{remark}
Suppose $T\ra S$ is a map of schemes.  If $X$ is an $\SP$-linear scheme 
over $S$, then the base-change $X\times_S T$ is an
$\SP$-linear scheme over $T$.  This follows from two facts.  First, 
smooth, proper schemes over $S$ pull back to smooth, proper
schemes over $T$.  Second, if $Z\inc X$ is a closed inclusion of
schemes over $S$ with complementary open subscheme $U = X - Z$, then
$Z\times_S T\ra X\times_S T$ is also a closed inclusion of schemes
over $T$ with open complement $U\times_S T$.
\end{remark}

\begin{lemma}
\label{lem:P-linear-product}
The class of $\SP$-linear schemes over $S$ is closed under finite products
in the category of schemes over $S$.
\end{lemma}

Recall that products in the category of schemes over $S$ are constructed
as fiber products over $S$.

\begin{proof}
Let $P$ be any smooth, proper scheme.  First we will show that
$P \times Q$ is $\SP$-linear whenever $Q$ is $\SP$-linear.
Let $\cC$ be the class of schemes $Q$ such that
$P \times Q$ is $\SP$-linear.  We want to show that $\cC$ contains
all $\SP$-linear schemes.  Note that $\cC$ contains $\SP$ because
the product of two smooth, proper schemes is again smooth and proper.
Now consider a closed inclusion $Z \inc X$ of smooth schemes 
with open complement $U$.
Then $P \times Z \map P \times X$ is a closed inclusion with open
complement $P \times U$.  If any two of $Z$, $X$, and $U$ belong to
$\cC$, then two of $P \times Z$, $P \times X$, and $P \times U$
are $\SP$-linear.  By definition, the third is also $\SP$-linear.
This shows that $\cC$ contains the class of $\SP$-linear schemes.

Now let $Q$ be a fixed $\SP$-linear scheme, and let $\cD$ be the class
of schemes $W$ such that $W \times Q$ is $\SP$-linear.  If we can show
that $\cD$ contains all $\SP$-linear schemes, then we will be done.
We have already shown that $\cD$ contains $\SP$.

Let $Z \map X$ be a closed inclusion of smooth schemes 
with open complement $U$.
Then $Z \times Q \map X \times Q$ is a closed inclusion of smooth schemes
with open complement $U \times Q$.  
The argument from the first paragraph of the proof applies again,
and $\cD$ contains all $\SP$-linear schemes.
\end{proof}

\begin{ex}
\label{ex:DQ}
If $R$ is a commutative ring in which $2$ is invertible, then the
scheme $DQ_n$ is $\SP$-linear over $\Spec R$.  It is the open
complement of the closed inclusion $Q_{n-1} \inc \P^n$, and both
$Q_{n-1}$ and $\P^n$ are smooth and proper.  It follows that the
scheme $DQ_m \times DQ_n$ is $\SP$-linear by Lemma
\ref{lem:P-linear-product}.
\end{ex}

Now let $V$ be a fixed strict Hensel local domain.  This means that
$V$ is a local integral domain which satisfies Hensel's lemma and
whose residue field $k$ is separably closed.  We also let $F$ be a
separably closed field containing the field of fractions of $V$.

\begin{ex}
\label{ex:Witt-vectors}
In our applications we will usually have
$k = \brF_p$, the algebraic closure of the finite field $\F_p$.
Also, we will have $V=W(\brF_p)$, the ring of Witt
vectors as defined in \cite[Sec.~II.6]{Se}.
Then $V$ is a
strict Hensel local domain whose residue field is $\brF_p$
by \cite[Thms.~II.5.3, II.5.5, II.6.8]{Se}.
Note that $V$ and $F$ have characteristic zero.
\end{ex}

The following theorem generalizes well-known properties of
smooth, proper schemes to $\SP$-linear schemes.  It will be applied
later in Proposition \ref{pr:compare}
to $DQ_n$ and $DQ_m \times DQ_n$ in order to analyze their
\'etale realizations.

\begin{thm}
\label{thm:P-linear-cohomology}
Suppose that $2$ is relatively prime to $\chara(k)$,
and let $X$ be $\SP$-linear over $\spec V$.
Then both maps in the diagram
\[
X_k \map X \leftarrow X_F
\]
induce isomorphisms in \'etale cohomology with coefficients in $\Z/2$.
\end{thm}

\begin{proof}
Let $\cC$ be the class of schemes $X$ over $\spec V$ such that
both maps induce \'etale cohomology isomorphisms.
We want to show that $\cC$ contains all $\SP$-linear schemes.

If $X$ is proper over $\Spec V$, then $H_{\et}^*(X;\Z/2) \ra
H_{\et}^*(X_k;\Z/2)$ is an isomorphism by \cite[Cor.~VI.2.7]{M}
or \cite[Cor.~XII.5.5]{SGA4}.
If $X$ is both smooth and proper then $H_{\et}^*(X;\Z/2) \ra
H_{\et}^*(X_F;\Z/2)$ is an isomorphism by 
\cite[proof of Cor.~VI.4.2]{M} or \cite[Cor.~XVI.2.2]{SGA4}.
Therefore $X$ belongs to $\cC$, and $\cC$ contains $\SP$.

Now suppose that $Z \inc X$ is a closed inclusion of smooth $V$-schemes
with open complement $U$.
We have a diagram
\[
\xymatrixcolsep{1.4pc}
\xymatrix{
\cdots & H_{\et}^p(U_F;\Z/2) \ar[l] & H_{\et}^p(X_F;\Z/2) \ar[l] & 
   H_{\et}^{p-2c}(Z_F; \Z/2) \ar[l] & \cdots \ar[l] \\
\cdots & H_{\et}^p(U;\Z/2) \ar[l]\ar[d]\ar[u] & 
   H_{\et}^p(X;\Z/2) \ar[l]\ar[d]\ar[u] & 
   H_{\et}^{p-2c}(Z; \Z/2) \ar[l]\ar[d]\ar[u] & \cdots \ar[l] \\
\cdots & H_{\et}^p(U_k;\Z/2) \ar[l] & H_{\et}^p(X_k;\Z/2) \ar[l] & 
   H_{\et}^{p-2c}(Z_k; \Z/2)) \ar[l] & \cdots \ar[l]     }
\]
in which the rows are long exact Gysin sequences (for the middle row,
see the remark following this proof).  Here $c$ is
the codimension of $Z$ in $X$.
If any two of $Z$, $X$, and $U$ belong to $\cC$, then both maps
in two of the columns
are isomorphisms.  By the five lemma, both maps in 
the third column are also isomorphisms,
so the third scheme belongs to $\cC$.  This shows that $\cC$ contains
all $\SP$-linear schemes.
\end{proof}

\begin{remark}
Beware that the middle row in the above diagram is a Gysin sequence over
a base scheme $\Spec V$ in which $V$ is not a field. 
By \cite[Thm.~VI.5.1]{M} (or \cite[XVI.3.7]{SGA4}), 
there is a long exact sequence 
\[
\cdots \leftarrow H_{\et}^p(U;\Z/2) \leftarrow H_{\et}^p(X;\Z/2) \leftarrow
H_{\et}^{p-2c}(Z; \cF) \leftarrow \cdots,
\]
where $\cF$ is a sheaf on $Z$ that is locally isomorphic to 
the constant sheaf $\Z/2$.
Since $\Aut(\Z/2)$ is trivial, it follows that $\cF$
must be the constant sheaf $\Z/2$.

This issue explains why the statement of 
Theorem \ref{thm:P-linear-cohomology} is specialized
to coefficients in $\Z/2$, rather than the usual twisted coefficients
$\mu_l(q)$ of $l$th roots of unity with $l$ relatively prime to $\chara(k)$.
\end{remark}

\begin{ex}
\label{ex:H-DQ}
Consider the scheme $DQ_n$ defined over $\brF_p$, the algebraic
closure of the finite field $\F_p$ with $p$ a prime greater than $2$.  
Recall the ring of Witt vectors $V = W(\brF_p)$ from
Example \ref{ex:Witt-vectors}.  Note that 
$\DQ_{n}$ is
defined over $V$ since it is in fact defined over $\Z$.

Let $F$ be the separable closure of the field of fractions of $V$.
Theorem \ref{thm:P-linear-cohomology} implies that 
$H_{\et}^*(DQ_n;\Z/2)$ is isomorphic to $H_{\et}^*((DQ_n)_F;\Z/2)$.
But $F$ has characteristic zero, so 
\cite[Thm.~III.3.12]{M} together 
with \cite[Cor.~VI.4.3]{M} 
(or \cite[Thm.~XI.4.4]{SGA4} together with
\cite[Cor.~XVI.1.6]{SGA4})
tells us that
$H_{\et}^*((DQ_n)_F;\Z/2)$ is isomorphic to the singular
cohomology ring $H^*(DQ_n(\C);\Z/2)$.  
Since $DQ_n(\C)$ is homotopy equivalent to $\RP^n$ by \cite[Lem.~6.3]{Lw}, 
we have computed that 
\[
H_{\et}^*(DQ_n;\Z/2) \cong \Z/2[x]/(x^{n+1}),
\]
where $x$ is a class of degree 1.

The same computation works for $\DQ_{n/k}$ where $k$ is any separably
closed field of characteristic different from 2.
Use \cite[Cor.~VI.4.3]{M} or \cite[Cor.~XVI.1.6]{SGA4} 
to reduce to the cases $k=\brF_p$ or $F = \C$,
depending on whether the characteristic is positive or zero, similarly
to the proof of Corollary~\ref{cor:2-wk-eq} below.
\end{ex}

\begin{ex}
\label{ex:H-DQxDQ}
Consider the scheme $DQ_m \times DQ_n$
over a separably closed field $k$ 
such that $\chara(k) \neq 2$.
As in Example \ref{ex:H-DQ}, 
we compute that $H_{\et}^*(DQ_m \times DQ_n;\Z/2)$ is isomorphic to
\[
H^*(\RP^m \times \RP^n;\Z/2) \cong \Z/2[x_1,x_2]/(x_1^{m+1},x_2^{n+1}),
\]
where $x_1$ and $x_2$ are classes of degree 1.
\end{ex}


\section{\'Etale realizations}
\label{sctn:et}

Recall that to any scheme $X$ one can associate a pro-simplicial set
$\Et(X)$ called the \dfn{\'etale realization} of $X$.  This is defined
in \cite[Def.~4.4]{F}, 
where it is called the `\'etale topological type'.
The construction of $\Et(X)$ is technical, but fortunately
the details are not necessary for the arguments below.

Note that if $X$ is a scheme over $\Spec R$ and $R\ra T$ is a ring map,
then the base-change map $X_T\ra X$ induces a map $\Et(X_T)\ra \Et(X)$.

Here are the facts that we will need about \'etale realizations.  They
are analogues of the results about \'etale cohomology used in
Section \ref{sctn:et-cohlgy}.

\begin{enumerate}[(1)]
\item \cite[Prop.~5.9]{F}: 
If $X$ is any scheme, then $H_{\et}^*(X;M)$ is naturally isomorphic
to $H^*(\Et(X);M)$ for any constant coefficients $M$.
In fact, a stronger statement concerning twisted coefficients
can be made, but we won't need it.
\item \cite[Thm.~8.4]{F} or \cite[Thm.~12.9]{AM}:
If $X$ is a scheme of finite type over $\C$, then there is a canonical
zig-zag  in $\pro\sSet$ of the form
\[ \Et(X) \la \SEt(X) \ra cX(\C)
\]
where both maps are $\Z/l$-cohomological equivalences, for any prime $l$.
Here $\SEt(\blank)$ is a certain functor defined in
\cite[Thm.~8.4]{F}.
\item
\cite[Prop.~8.6, Prop.~8.7]{F} or
\cite[Cor.~12.13]{AM}:
Suppose that $V$ is a strict Hensel local domain---that is, a local domain
satisfying Hensel's lemma whose residue field $k$ is separably
closed.  Let $F$ be a separably closed field containing the quotient
field of $V$.  If $X$ is a scheme that is both smooth and proper over 
$\Spec V$, then the maps
of pro-spaces
\[ \Et(X_k) \ra \Et(X) \la \Et(X_F)
\]
are both $\Z/l$-cohomological equivalences, for any prime $l$
different from $\chara(k)$.  
\end{enumerate}

We need one additional property of \'etale realizations, concerning
their behavior with respect to field extensions.  

\begin{lemma}
\label{lem:et-sep-closed}
Let $E\inc F$ be an inclusion of separably closed fields.  If
$X$ is a scheme over $E$, then the canonical map $\Et(X_F) \ra
\Et(X)$ is a $\Z/l$-cohomological equivalence for any prime $l$
different from $\chara(F)$.  
\end{lemma}

\begin{proof}
By the definition of $\Z/l$-cohomological equivalences and property (1),
we need only observe that the map
$H_{\et}^*(X;\Z/l) \map H_{\et}^*(X_F;\Z/l)$
is an isomorphism by \cite[Cor.~VI.4.3]{M} or \cite[Cor.~XVI.1.6]{SGA4}.
\end{proof}

The careful reader will notice that this lemma is slightly
different than \cite[Cor.~12.12]{AM}.
We have not assumed that $X$ is proper; instead, we have
assumed that $l$ is relatively prime to the characteristics of the fields.

We wish to
compute $\Et(DQ_n)$ and $\Et(DQ_m\times DQ_n)$ and to show that
they can be connected by $\Z/2$-cohomological equivalences to the
constant pro-spaces $c\RP^n$ and $c(\RP^m \times \RP^n)$.  This will
go a long way towards establishing Proposition~\ref{pr:mainaxialp}.  
As the
scheme $DQ_n$ can be lifted to characteristic zero, and we know
that $DQ_n(\C)$ is homotopy equivalent to $\RP^n$ by \cite[Lem.~6.3]{Lw},
properties (1)--(3)
and Lemma \ref{lem:et-sep-closed}
can almost be used to provide the required equivalences.
Although $DQ_n$ is smooth, it is not proper; thus property
(3) does not apply directly.  However, since
$DQ_n$ is $\SP$-linear, the results of Section \ref{sctn:et-cohlgy}
will get around this problem.

Now assume that there is a diagram 
\begin{myequation}
\label{eq:rings}
\xymatrix{
A \ar[r]\ar[d] & V \ar@{ >->}[r] & E \ar@{ >->}[d]\\
\C \ar@{ >->}[rr] && K
}
\end{myequation}
of rings in which the indicated maps are inclusions and where 
\begin{itemize}
\item $V$ is a strict Hensel local domain with separably closed residue
field $k$,
\item $E$ and $K$ are separably closed fields
(and so $E$ is an extension of the quotient field of $V$).
\end{itemize}
If $X$ is a scheme over $A$, one obtains schemes 
$X_V$, $X_E$, $X_{\C}$, $X_K$, and $X_k$ via base change.
Let $X_\C(\C)$ denote the topological space of $\C$-valued points of
$X_\C$, equipped with the analytic topology.

\begin{prop}
\label{pr:compare}
Suppose given a diagram as in (\ref{eq:rings}) such
that the residue field $k$ of $V$ has characteristic different from $2$.
For any scheme $X$ over $A$, there is a natural diagram
\[
\Et(X_{k}) 
\ra \Et(X_V) \la \Et(X_E) \la \Et(X_K) \ra \Et(X_\C) \la \SEt(X_\C) \ra 
cX_{\C}(\C)
\]
of pro-spaces. If $X$ is $\SP$-linear, then
each map is a $\Z/2$-cohomological
weak equivalence.  Here $\SEt$ is a functor from schemes over $\C$
to pro-spaces defined in 
\cite[Thm.~8.4]{F}.  
\end{prop}

Recall the $\Z/2$-cohomological homotopy category $\ho_{\Z/2}(\pro \sSet)$
of pro-spaces that is briefly described in Section \ref{sec:main}
and is discussed more carefully in Section \ref{sctn:pro}.
The point of the proposition is that when $X$ is $\SP$-linear there is
a natural isomorphism between $cX_{\C}(\C)$ and $\Et(X_k)$ in
$\ho_{\Z/2}(\pro \sSet)$.

\begin{proof}
The first two maps are $\Z/2$-cohomological weak equivalences as
a result of Theorem \ref{thm:P-linear-cohomology}
and property (1).
The middle two maps $\Et(X_K) \ra \Et(X_E)$ and $\Et(X_K)\ra
\Et(X_\C)$ are $\Z/2$-cohomological equivalences by 
Lemma \ref{lem:et-sep-closed}.
Finally, the last two maps are $\Z/2$-cohomological weak equivalences 
because of \cite[Thm.~8.4]{F}.
\end{proof}

\begin{cor}
\label{cor:2-wk-eq}
Let $k$ be a separably closed field of characteristic not equal to 2,
and consider the schemes $\DQ_n$ and $DQ_m \times DQ_n$
defined over $k$.
In $\ho_{\Z/2}(\pro \sSet)$,
there are isomorphisms
\[
c\RP^{n} \map \Et(DQ_{n})
\]
and
\[
c(\RP^{m} \times \RP^{n}) \map \Et(DQ_{m} \times DQ_{n}).
\]
\end{cor}

\begin{proof}
First suppose that $k$ has characteristic zero.
Choose an embedding $\bar{\Q}\inc k$, and note that $DQ_{n/k}$ is the base
extension of $\DQ_{n/\bar{\Q}}$.  There is a zig-zag of maps
\[ \Et(DQ_{n/\C}) \ra \Et(DQ_{n/\bar{\Q}}) \la \Et(DQ_{n/k}) \]
which are both $\Z/2$-cohomological equivalences by
Lemma~\ref{lem:et-sep-closed}.  Finally, property (2) from the
beginning of this section says that $\Et(DQ_{n/\C})$ is
$\Z/2$-cohomologically equivalent to $c(\RP^n)$, since
$\DQ_{n}(\C)$ is homotopy equivalent to $\RP^n$ by \cite[Lem.~6.3]{Lw}.
The same argument works for $DQ_m\times DQ_n$.

Now suppose that $k$ has characteristic $p$, where $p > 2$.
Choose an embedding
${\F}_p^s\inc k$ where $\F_p^{s}$ 
is the separable closure of
the finite field $\F_p$.  
Note that $DQ_{n/k}$ is the base extension of
$DQ_{n/{\F_p^{s}}}$, so  Lemma~\ref{lem:et-sep-closed} gives us a
$\Z/2$-cohomological equivalence
\[ \Et(DQ_{n/k}) \ra \Et(DQ_{n/\F^s_p}). \]
It therefore suffices to consider the case $k=\F^{s}_p$.

Next consider the inclusion $\F^{s}_p \inc \brF_p$, where
$\brF_p$ is the algebraic closure of $\F_p$.  This is an inclusion
of separably closed fields, so
Lemma~\ref{lem:et-sep-closed} again gives us a
$\Z/2$-cohomological equivalence
\[ \Et(DQ_{n/\brF_p}) \ra \Et(DQ_{n/\F_p^s}). \]
So it further suffices to consider the case $k=\brF_p$.

Recall from Example \ref{ex:Witt-vectors} that
$V = W(\brF_p)$ is the strict Hensel local domain
of Witt vectors for $\brF_p$.
Let $A = \Z$,
and let $E$ be a separable closure of the quotient
field of $V$.  For $K$ we choose any field containing both $E$ and
$\C$, such as the compositum of $E$ and $\C$.

We are now in a position to apply the preceding proposition, since
$\DQ_n$ is $\SP$-linear.  We find that $\Et(\DQ_{n/\brF_p})$
is isomorphic to 
$c\DQ_{n/\C}(\C)$ in the $\Z/2$-cohomological homotopy category.  
We have already observed many times that
the spaces $\DQ_{n/\C}(\C)$ and $\RP^n$ are
homotopy equivalent, so this completes the proof for $\DQ_n$.

Everything works the same for the product $\DQ_m\times \DQ_n$, using
that this scheme is $\SP$-linear by Lemma~\ref{lem:P-linear-product}.
\end{proof}

\begin{remark}
\label{re:commsquare}
We will need the following fact later in the proof
of Proposition \ref{pr:mainaxialp} given at the end of this section.
Consider the rational point $*=[1:0:0:\cdots:0]$ of 
the scheme $DQ_n$ defined over a field $k$.  We regard
$*$ as a map of schemes $\Spec k \ra DQ_n$.
There is an induced map $j_1= \id \times *\colon DQ_m \ra DQ_m\times DQ_n$.
We then have a diagram in $\Ho_{\Z/2}(\pro
\sSet)$  of the form
\[\xymatrix{
c(\RP^m) \ar[r]^\iso \ar[d] & \Et(DQ_{m}) \ar[d]^{\Et(j_1)} \\
c(\RP^{m}\times \RP^n) \ar[r]^-\iso & \Et(DQ_m\times DQ_n)}
\]
where the two isomorphisms are from Corollary \ref{cor:2-wk-eq} 
and the left
vertical map comes from the evident inclusion $\RP^m \inc \RP^m\times\RP^n$
(using the point $[1:0: \cdots:0]$ of $\RP^n$).
This diagram is commutative because of the
naturality of the maps in Proposition~\ref{pr:compare}.

The same observations apply to the evident map
$j_2\colon DQ_n \map DQ_m \times DQ_n$.
\end{remark}

We are now ready to complete the last remaining detail of the paper.

\begin{proof}[Proof of Proposition~\ref{pr:mainaxialp}]
Assume that a sums-of-squares formula of type $[r,s,n]$ exists over a
field $k$ with $\chara(k)\neq 2$.
Let $k^s$ be the separable closure of $k$, and note that the sums-of-squares
formula over $k$ automatically exists over $k^s$ also. 
 Hence we may assume that $k$ is separably closed.  

As explained in \cite[Sec.~1]{DI1}, the sums-of-squares formula
induces a map $f\colon DQ_{r-1}\times DQ_{s-1} \ra DQ_{n-1}$.  At the
level of points, this map sends the pair
$([x_1:\ldots:x_r],[y_1:\ldots:y_s])$ to $[z_1:\ldots:z_n]$, using the
notation from Section~\ref{sctn:intro}.

Let $*=[1\colon 0 \colon 0 \colon \cdots \colon 0]$ 
be the basepoint in both $DQ_{r-1}$
and $DQ_{s-1}$.  Let $j_1\colon DQ_{r-1}\inc DQ_{r-1}\times DQ_{s-1}$
be the closed inclusion as defined in Remark \ref{re:commsquare}.
Define $j_2$ similarly.  

Recall from Example~\ref{ex:H-DQ} that the group $H_{\et}^1(DQ_{n-1};
\Z/2)$ is isomorphic to $\Z/2$; let $x$ be the generator.
Let $x_1$ and $x_2$ be the classes
$\pi_1^*(x)$ and $\pi_2^*(x)$ of
$H_{\et}^1(DQ_{r-1} \times DQ_{s-1};\Z/2)$,
where $\pi_1$ and $\pi_2$ are the
obvious projections.  Note that $j_1^*(x_1)=x = j_2^*(x_2)$
and $j_2^*(x_1)=0 = j_1^*(x_2)$.
This shows that $x_1$ and $x_2$ are distinct
non-zero classes.  But by Example~\ref{ex:H-DQxDQ}, $H_{\et}^1(DQ_{r-1}
\times DQ_{s-1};\Z/2)$ is isomorphic to $\Z/2 \times \Z/2$, so $x_1$
and $x_2$ are generators for this group.

Since the sums-of-squares formula is bilinear, 
the compositions $fj_1\colon DQ_{r-1} \map DQ_{n-1}$ and
$fj_2\colon DQ_{s_1} \map DQ_{n-1}$ are linear inclusions.
This shows that $f^*(x) = x_1 + x_2$.
(compare \cite[Prop.~2.5]{DI1}, which is the analogous result for
motivic cohomology).

According to Lemma \ref{lem:rep-H1},
for any scheme $X$, the set 
$[\Et(X),c\RP^\infty]_{\Z/2}^{\pros}$ of homotopy classes 
in $\ho_{\Z/2}(\pro \sSet)$
is naturally
isomorphic to $H^1(\Et(X);\Z/2)$,
which in turn is isomorphic to $H_{\et}^1(X;\Z/2)$ 
by property (1) from the beginning of Section \ref{sctn:et}.
The element $x$ of 
$H_{\et}^1(DQ_{n-1};\Z/2)$ therefore
gives us a homotopy class
$\Et(DQ_{n-1}) \map c\RP^\infty$ in $\ho_{\Z/2}(\pro \sSet)$.
Similarly, the elements 
$x_1$ and $x_2$ of
$H_{\et}^1(DQ_{r-1} \times DQ_{s-1}; \Z/2)$ give us
a homotopy class 
$\Et(DQ_{r-1} \times DQ_{s-1}) \map c(\RP^\infty\times \RP^\infty)$
in $\ho_{\Z/2}(\pro \sSet)$; here we are using the isomorphism between
\[ H^1(DQ_{r-1}\times DQ_{s-1};\Z/2) \times H^1(DQ_{r-1}\times
DQ_{s-1};\Z/2) \]
and
\[ [\Et(DQ_{r-1}\times
DQ_{s-1}),c(\RP^\infty\times\RP^\infty)]^{\pros}_{\Z/2}
\]
from the
second part
of Lemma \ref{lem:rep-H1}.

We now have a diagram
\[ \xymatrix{ \Et(DQ_{r-1}\times DQ_{s-1}) \ar[r] \ar[d]_{\Et(f)} &
c(\RP^\infty \times \RP^\infty) \ar[d]^{c(\mu)} \\
\Et(DQ_{n-1}) \ar[r] & c\RP^\infty
}
\]
in $\ho_{\Z/2}(\pro\sSet)$.  
On cohomology, the left vertical map takes $x$ to $x_1 + x_2$.
Commutativity of this diagram follows
from this and from the fact that $\mu$
represents the element $y_1 + y_2$
of $H^1(\RP^\infty \times \RP^\infty;\Z/2)$,
where $y$ is the generator of $H^1(\RP^\infty;\Z/2)$
and $y_1$ and $y_2$ are the classes $\pi_1^*(y)$ and $\pi_2^*(y)$
in $H^1(\RP^\infty\times\RP^\infty;\Z/2)$.

Using the maps supplied by Corollary \ref{cor:2-wk-eq},
we obtain a commutative diagram
\[ \xymatrix{ c(\RP^{r-1} \times \RP^{s-1}) \ar[r]^-\sim\ar[d] & 
   \Et(DQ_{r-1}\times DQ_{s-1}) \ar[r] \ar[d] &
c(\RP^\infty \times \RP^\infty) \ar[d]^\mu \\
c\RP^{n-1}\ar[r]^\sim & \Et(DQ_{n-1}) \ar[r] & c\RP^\infty
}
\]
in $\ho_{\Z/2}(\pro\sSet)$, where the left vertical map is defined using
the homotopy-inverse of the map $c\RP^{n-1} \map \Et(DQ_{n-1})$.
The outer square in this diagram is the desired commutative square.
We still need to show that the horizontal
maps correspond to the canonical inclusions; this identification
is necessary for the arguments in Section~\ref{sec:main}.

For the bottom composite,
there are only two homotopy classes
$c\RP^{n-1} \map c\RP^\infty$, and they are classified by the map
that they induce on $H^1(\blank;\Z/2)$ (using Lemma~\ref{lem:rep-H1}).
We have chosen $\Et(DQ_{n-1}) \map c\RP^\infty$ so that the
bottom composite is the non-zero homotopy class.  On the other hand,
any linear inclusion $\RP^{n-1} \map \RP^\infty$ also represents 
the non-zero homotopy class.

For the top composite, there are exactly sixteen homotopy
classes of maps $h\colon c(\RP^{r-1}\times\RP^{s-1}) \ra c(\RP^\infty\times
\RP^\infty)$, and again they are classified by the 
induced map on
$H^1(\blank;\Z/2)$.  More concretely, 
such maps are classified by the
values of the four expressions
\[
j_1^* h^* \pi_1^*(y), \quad
j_1^* h^* \pi_2^*(y), \quad
j_2^* h^* \pi_1^*(y), \quad \text{and} \quad
j_2^* h^* \pi_2^*(y) \quad
\]
(each of which has exactly two possible values).
Here $j_1$, $j_2$, $\pi_1$, and $\pi_2$ are the evident 
inclusions and projections, as above.

Let $h$ denote the top composite in the above diagram. Check
that $j_1^*h^*\pi_1^*(y)$ 
and $j_2^*h^*\pi_2^*(y)$ are non-zero, while
$j_1^*h^*\pi_2^*(y)$ and $j_2^*h^*\pi_1^*(y)$ are zero;
this step uses the commutativity of the square in
Remark~\ref{re:commsquare}.  
Finally, if 
\[ j\times j\colon \RP^{r-1}\times \RP^{s-1}\ra \RP^\infty
\times \RP^\infty
\]
 is the standard inclusion, observe that 
$j_1^*(j\times j)^*\pi_1^*(y)$ and
$j_2^*(j\times j)^*\pi_2^*(y)$ are non-zero, while
$j_1^*(j\times j)^*\pi_2^*(y)$ and
$j_2^*(j\times j)^*\pi_1^*(y)$ are zero.
By the observation in the previous paragraph,
it follows that $h$ and $j \times j$ represent the same map in
$\ho_{\Z/2}(\pro \sSet)$.  This completes the proof.
\end{proof}

\begin{remark}
There is another, more elegant method for proving 
Proposition~\ref{pr:mainaxialp},
but we have avoided it because it requires more machinery.
The idea is to prove a motivic version of the proposition.  Let $k$
be a field, and recall that Morel and Voevodsky \cite{MV} have constructed
a model category representing the motivic homotopy theory
of smooth schemes over $k$.  

If $\chara(k)\neq 2$ and a sums-of-squares formula
of type $[r,s,n]$ exists over $k$, then
there exists a diagram 
\[
\xymatrixcolsep{2.4pc}
\xymatrix{ 
DQ_{r-1}\times DQ_{s-1} \ar[d]
   \ar[r]^-{j\times j}
& 
DQ_\infty \times DQ_\infty \ar[d]
\\
DQ_{n-1} \ar[r]^{j} 
& 
DQ_\infty }
\]
in the motivic homotopy category, where the horizontal maps are
standard inclusions and $DQ_{\infty}$ is the colimit (or ascending
union) of the schemes $DQ_n$ for all $n$.

The paper \cite{I2} constructs a functor $\Et$
from the motivic homotopy category to 
$\ho_{\Z/2}(\pro\sSet)$.
This functor extends the usual \'etale realization on schemes,
and applying it to the above diagram yields a diagram of pro-spaces.
Some work is required to identify $\Et(DQ_\infty)$ with $c\RP^\infty$ and
$\Et(DQ_\infty\times DQ_\infty)$ with $c(\RP^\infty\times \RP^\infty)$ (and
the map between them with $\mu$), but this then yields the desired diagram of
Proposition~\ref{pr:mainaxialp}.  

Another relevant tool for this approach is
the stable \'etale realization of \cite{Qk} from
the motivic stable homotopy category to another stable homotopy category
that is closely related to the $\Z/2$-cohomological homotopy category
of pro-spectra.  
\end{remark}


\bibliographystyle{amsalpha}

\end{document}